\documentstyle[12pt]{article}%
\def\nc{\newcommand}
\def\rnc{\renewcommand}
\nc{\nn}{\nonumber}
\makeatletter
\@addtoreset{equation}{section}
\rnc{\theequation}{\arabic{section}.\arabic{equation}}
\makeatother
\nc{\ignore}[1]{}
\nc{\be}{\begin{equation}}
\nc{\ee}{\end{equation}}
\nc{\ba}{\begin{array}}
\nc{\ea}{\end{array}}
\nc{\bea}{\begin{eqnarray}}
\nc{\eea}{\end{eqnarray}}

\rnc{\a}{\alpha}
\rnc{\b}{\beta}
\rnc{\d}{\delta}
\nc{\D}{\Delta}
\nc{\Db}{\bar\D}
\nc{\e}{\epsilon}
\nc{\eb}{\bar{\epsilon}}
\nc{\f}{\phi}
\nc{\fb}{\bar{\phi}}
\nc{\vf}{\varphi}
\nc{\p}{\psi}
\rnc{\c}{\chi}
\nc{\cb}{\bar{\c}}
\nc{\Qb}{\bar{Q}}
\rnc{\l}{\lambda}
\rnc{\o}{\omega}
\rnc{\t}{\theta}
\nc{\T}{\Theta}
\nc{\tb}{\bar{\theta}}
\nc{\im}{\imath}
\nc{\s}{\sigma}
\nc{\gb}{\bar g}
\def\u5{{\cal U}_5}
\def\Qb{{\bar Q}}
\def\Qb{{\bar Q}}

\nc{\trac}[2]{{\textstyle\frac{#1}{#2}}}
\nc{\half}{\trac{1}{2}}

\nc{\mat}[4]{\left(\begin{array}{cc}#1&#2\\#3&#4\end{array}\right)}

\def\Tr{\mathop{\rm Tr}\nolimits}
\nc{\ot}{\otimes}
\rnc{\ss}{\subset}
\nc{\ra}{\,\rangle}
\nc{\la}{\langle\,}
\nc{\ad}{a^\dagger}
\nc{\cd}{c^\dagger}
\nc{\rf}[1]{(\ref{#1})}
\nc{\tbm}{{\tilde B}_-}
\nc{\tbp}{{\tilde B}_+}
\nc{\tbpi}{{\tilde B}_{+i}}
\nc{\tbmi}{{\tilde B}_{-i}}
\nc{\tbpj}{{\tilde B}_{+j}}
\nc{\tbmj}{{\tilde B}_{-j}}
\nc{\del}{\partial}
\rnc{\=}{ \; = \; }
\rnc{\-}{ \: - \: }
\rnc{\+}{ \: + \: }
\nc{\qq}{\qquad}
\nc{\eq}[1]{(\ref{#1})}
\nc{\com}[3]{[ \: #1 \: ,\: #2 \: ] \= #3}
\nc{\I}{${\tt I}${ }}
\nc{\II}{${\tt I\!I}${ }}
\nc{\III}{${\tt I}^\prime${ }}
\nc{\IV}{${\tt I\!I}^\prime${ }}
\nc{\mr}[1]{(\ref{#1})}
\nc{\ie}{{\em i.e.}{}}
\begin{document}

\begin{titlepage}
\title{{\bf Equivalence of $q$-bosons \\
using the exponential phase operator }}
\date{}
\author{
{\Large S. U.  Park}\thanks{Permanent address:
Department of Physics, Jeonju University,
1200 Hyoja-3 Chonju, Chonbuk, 560-759, Korea}
\\[1cm]
Department of Mathematics \\
School of Mathematical Sciences\\
Australina National University\\
Canberra  ACT 0200\\
Australia
}
\maketitle

\vspace{1cm}
\begin{abstract}
Various forms of the $q$-boson are explained and their hidden symmetry revealed
by transformations using
the exponential phase operator.
Both the one-component and the multicomponent $q$-bosons are discussed.
As a byproduct, we obtain a new boson algebra having a shifted vacuum structure
and define a global operatal $U(1)$ gauge transformation.
\end{abstract}

\end{titlepage}

\section{Introduction}
The $q$-boson\footnote{There are equivalent terminologies such as
the $q$-Heisenberg-Weyl algebra and the $q$-oscillator.
But we take the $q$-boson
in order to emphasize its relation with the complex variable.}
is the simplest $q$-deformed algebra having a deformation
parameter. Various kinds of $q$-bosons have been introduced for their own
motives. They can be transformed into each other by proper redefinitions of
their generators. The existence of such transformations suggests that
the $q$-bosons may be equivalent to each other,
and there will exist underlying symmetries that make the equivalence
sensible.
However, we cannot see any clues for finding such symmetries.

One of the reasons for this fallacy of understanding may be connected to
the difference in the methods of construction, such as
in the process of the Schwinger realization of $su_q(2)$ in [1,2,3],
and of $su_q(n)$ in \cite{sun},
as the components of the fundamental representation of $su_q(n)$ \cite{pusz},
or for other motives [5 - 10].
Each type of $q$-bosons transform into each other [11,12,13].

Another reason is that the Hopf algebraic structure does not fix the $q$-boson
uniquely [14,15]. However, in this paper we will not be concerned about
the Hopf structure [16,17].

Recently, a new method of constructing $q$-bosons has
been presented in [18, 19] in which $q$-bosons are rewritten by the biproduct
of different generators of bosonic type and taking an expectation value for
a density operator. Although its connection with
the Hopf structure has not been clarified, it explains some $q$-properties easily.

In this paper, we restrict the method to the boson and
the exponential phase operator [20, 21, 22].
Depending on the explicit choice, the different forms of $q$-bosons appear,
clearly pointing to the difference and similarity between the $q$-bosons.
As a result, $q$-bosons, treated in this paper, are related to each other
with the operator version of $U(1)$ gauge transformation
and the similarity transformation by the exponential phase operator.
As a by-product, we find a new boson algebra which is equivalent to
the normal boson algebra except of the shifted vacuum,
and the operator version of the $U(1)$ gauge transformation
which becomes a normal gauge transformation in the classical limit.
We work both on the one-component $q$-bosons and on the multicomponent forms.

The paper is arranged as follows. In the next section, we treat the statistical mixture
in order to see the physical meaning of the exponential phase operator
and to define the physical density operator for $q$-bosons.
In the third section, we construct $q$-bosons and find connections between them.
We also define a new boson algebra and the operator version of the $U(1)$ gauge
transformation. In the fourth section, we treat the multi-component $q$-bosons,
and find a hidden symmetry. Finally, we briefly discuss problems along the way,
encountered along with extensions to the $q$-fermion and so on.

\section{The statistical mixture and the Cuntz algebra}

Let us consider the statistical mixture of the bosonic system in order
to define the density operator of the expectation value which will
appear in the process of constructing $q$-bosons.
We also discuss the physical meaning of the exponential phase operator.
A more detailed physical system related with this section should be found
in  quantum optics [20].

Assume that a system has a known probability $P_R$ when it is in the state
$\mid R\ra$.
Here, $R$ is a label that runs over a set of pure states sufficient to describe
the system. The states described by the probability $P_R$  are called  the
{\em statistical mixture}, and the magnitude of the $P_R$ for $\mid R \ra$
contains all the available information about the state.
An example of a statistical mixture is provided by the thermal excitation of
the photons in a cavity mode.
The probability $P_n$ means that $n$ photons are excited at the temperature $T$.
In this case, the result of an experiment depends on an ensemble average of
some observable quantity.

Consider some quantum mechanical operator $O$. The average value of
the observable for the pure state $\mid R \ra$ is $\la R\mid O \mid R \ra$,
and hence the ensemble average of the observable for the statistical mixture
specified by $P_R$ is
\be
\la O \ra \= \sum_R P_R \la R \mid O \mid R\ra.
\ee
It will be assumed that $P_R$ is a normalized probability distribution
\be
\sum_R P_R \= 1.
\ee
The average $\la O \ra$ is independent of the particular complete set of
the state chosen for evaluation. This fact is apparent by defining the density
operator $\rho$ as follows
\be
\rho\=\sum_R P_R \mid R\ra \la R \mid.
\ee
The density operator contains exactly the same information as the probability
distribution $P_R$. The average $\la O \ra$ can be written
\be
\la O \ra \= \Tr (\rho O ),
\ee
where the trace of an operator ( hence abbreviated to $\Tr$ ) is the sum of its
diagonal matrix elements for any complete set of states. In particular, the
expectation value for the identity operator is equal to  unity by the
definition of the probability
$ \Tr(\rho)\= \sum_R P_R \= 1.$

We can regard a pure state as a special case of a statistical mixture in which
one of the probabilities $P_R$ is equal to  unity, and all the remaining
$P_R$ are zero. The pure-state density operator is defined as
\be
\rho_{R} \= \mid R\ra \la R\mid.
\ee
In this case, the state is retracted on a particular state, and statistical
description becomes somewhat redundant. However, the concept of the density
operator remains valid. Also the pure state density operator satisfies the
property of the projection operator
\be
\rho^2\!\!{}_R \= \rho_R,
\ee
which is easily proved from the definition (2.5).

From now on, for later application, we will particularly restrict our attention
to the bosonic system $(a,\ad, N)$, where
\be
[ a, \ad ]\= 1, \qquad [N,\ad]\= \ad.
\ee
The Hilbert state $\mid n \ra$ is characterized by the number operator
\be
N\mid n\ra \= n \mid n\ra ,\qquad a\mid n \ra \= \sqrt{n} \mid n-1\ra .
\ee
The ground state $\mid 0 \ra$ is defined by the annihilation operator
$a \mid 0 \ra =0$, and all other states are constructed from it:
$\mid n\ra \= \frac{1}{\sqrt{n!}} (\ad)^n \mid 0\ra $.

The generators of the boson algebra can be decomposed into the magnitude $N$,
\ie the number operator and the phase $\f$,  if we quantize in the polar
coordinate $(N,\f)$ of the phase space,
\be
a\= e_a\:\sqrt{N}, \qq \ad\= \sqrt{N}\:e^\dagger_a.
\ee
Here, the newly introduced operators $(e_a, e_a^\dagger)$ are called the exponential
phase operators. To denote its relation with the $a$-boson, we write
subscript $a$. There are various different definitions of the exponential
phase operators, but we will choose the Susskind-Glogower notation [21, 22],
\be
e_a e^\dagger_a \= 1,\qq e^\dagger_a e_a \= 1 \- \mid 0 \ra \la 0 \mid,
\ee
which is not unitary.
The algebra is called the Cuntz algebra, and will act an important role in the
large $N$-expansion of the matrix model [23, 24, 25]
and to calculate the anomaly [26].
Although the Susskind-Glogower notation is mathematically easy and useful to treat,
it is uncertain that it is a real physical operator which describes the nature [27].
Neverthless we will choose the notation for its convenience.

As we see in products of the exponential phase operators, they
are similar to the real phase (the unitary operator) except the vacuum.
The commutation relations among $(N,e_a^\dagger, e_a)$ are derived from
the boson algebra (2.7) and the definition (2.9),
\be
[N, e_a] \= -e_a,\qq
[N, e_a^\dagger] \= e_a^\dagger.
\ee
The exponential phase operator and its ajoint shift up and down one step
in the number state,
\be
e_a^\dagger \mid n \ra \= \mid n+ 1\ra, \qq
e_a \mid n \ra \= \mid n-1\ra.
\ee

We now turn to the properties of the exponential phase operators
in order to find their physical meanings [20].
They are easily seen in the coherent state of the boson algebra defined
by the eigenstate of the annihilation operator,
\be
a \mid z \ra \= z \mid z \ra.
\ee
This $\mid z \ra$ is normalized, $\la z\mid z \ra \=1$,
and thus expressed by the superposition of the number states
\be
\mid z \ra \= e^{-|z|^2/2} \sum_{n=0}^\infty \frac{z^n}{\sqrt{n!}} \mid n \ra.
\ee
The coherent density operator $\rho_{c}$ can be written as
\be
\rho_{c} \= \mid z \ra \la z \mid.
\ee
The density operator shows a Poisson distribution with respect to the number
state.
\be
\mid \la n \mid z \ra \mid^2 \= \exp (-|z|^2) \frac{|z|^{2n}}{n!}.
\ee
The mean value of the number operator, (or mean number in short),
under this probability is $|z|^2$, \ie the radius
or the magnitude of the complex variable $z$,
\be
\la z \mid N \mid z \ra \= |z|^2.
\ee
Similarly, the expectation value of the exponential phase
operator becomes the physical (real) phase in the limit of the large mean number
( \ie the classical limit),
\be
\la z \mid e^{i\phi} \mid z \ra \= z \exp(-|z|^2)
\sum_n \frac{|z|^{2n}}{((n+1)!n!)^{1/2}}.
\ee
It is not possible to evaluate the summation analytically,  the asymptotic
expansion is obtained for large $|z|^2$, \ie the classical limit, with
\be
\sum_n \frac{|z|^{2n}}{n!(n+1)^{1/2}} \= \frac{\exp |z|^2}{|z|}
( 1\-\frac{1}{8|z|^2} \+ \cdots), \qquad |z|^2 \gg 1.
\ee
Thus, the expectation value of the exponential phase operator becomes
\be
\la z \mid e_a \mid z \ra
 \= \frac{z}{|z|}( 1\-\frac{1}{8|z|^2} \+ \cdots)\:\rightarrow\: e^{i\f}.
\ee
Since $z$ is a complex number, we can decompose it into its magnitude and phase,
$z\=|z|e^{i\f}$, then the exponential phase operator becomes the physical phase
asymptotically.

We define other well known density operators for later convenience.
The pure-state density operator for a number state is given by
\be
\rho_n \= \mid n \ra \la n \mid.
\ee
For the thermally excited state, we define the exponential distribution
(a Planck distribution ) for a given temperature $T$ and a quantal energy $\e_0$
as
\bea
\rho_e &=& (1-q^2) \sum_n q^{2n} \mid n \ra \la n \mid, \nn \\
       &=& (1-q^2) q^{2 N},
\eea
where a $q$-parameter is given by
\be
q^2\= \exp (-\e_0/k_B T).
\ee
This parameter will act like the $q$-parameter in a $q$-deformation.
Thus we can see a physical deformation depending on the temperature $T$ and
the quantal energy $\e_0$.

The following mean values under the exponential distribution (2.22)
are obtained as
\bea
\la \ad a \ra &=& \frac{q^2}{1\-q^2}, \\
\la a \ad \ra &=& \frac{1}{1\-q^2}.
\eea
Also those related with the exponential phase operators (2.10) are obtained as follows,
\bea
\la e_a^\a e_a^{\dagger\a} \ra &=& 1, \\
\la e_a^{\dagger\a} e_a^\a  \ra &=& q^{2\a}, \\
\la \t(N \- \a) \ra &=& q^{2\a}.
\eea
Here, $\t (x)$ is a step function and $\a\ge 0$ is taken as an integer
in order to keep the structure of the number state.

\section{Construction of $q$-bosons and the Hidden Transformations}

\nc{\bd}{b^\dagger}
We consider a system in which the algebra is described by the generators
$(D_\pm, D_0)$ satisfying the formal commutation relation
\be
[D_-, D_+] \= D_0.
\ee
We will not specify the commutators between $D_\pm$ and $D_0$. Although they
are important to close an algebra, their effect is null in the expectation value.
Furthermore, we assume that two sets of independent bosons
$\{ (a,\ad,N_a),(b,\bd,N_b)\}$  contribute to  the generators $D_\pm$
in the product forms.
The Schwinger representation of $su(1,1)$ ( or $su(2)$) is the most typical
example in which a product of two independent bosons forms the new algebra.
The simple product of two independent bosons
has been well studied [28],
so we should extend it into that of independent bosonic contributions
such as boson generators  themself, their exponential phase operators,
and combinations of the exponential phase operators and the number operators.
A part of the algebra has already been treated in [19].

We choose the forms of two operators $D_\pm$ as follows,
\be
D_\pm \= A_\pm B_\pm .
\ee
The $A_\pm$ and $B_\pm$ are assumed to be independent, \ie commuting,
and have the property of only $a$- and $b$-bosons respectively.
In order not to complicate matters further, we restrict the $A_\pm$
and $B_\pm$  to have only one quantal number of their own.
\bea
[ N_a, A_\pm ] \= \pm A_\pm, & [N_b, A_\pm] \= 0, \nn \\[1mm]
[N_a, B_\pm ] \= 0, & [N_b , B_\pm ] \= \pm B_\pm, \\[1mm]
[ A_\pm, B_\pm ] \=0, & [A_\mp, B_\pm] \= 0. \nn
\eea
Here $N_i,\; i=a,b$ is the number operator of each boson.

From the product form (3.2) of $D_\pm$, we rewrite the commutator (3.1) as follows,
\be
A_- A_+  B_- B_+ \-  A_+ A_-  B_+ B_-  \=  D_0 .
\ee
We assume that its full Hilbert space of $A_\pm$ in $D_\pm$ is known.
The property of the $B_\pm$ should be determined to keep the algebra
of $D_\pm$ depending on the forms of $D_0$.
Under this assumption, the equation (3.4) itself is the relation for $B_\pm$.
To solve it, we should take an expectation value for a density operator,
\be
\la A_- A_+ \ra \la B_- B_+ \ra \- \la A_+ A_- \ra \la B_+ B_- \ra
\= \la D_0 \ra.
\ee
We will discuss the results with respect to density operators.
Two types of density operators (2.21) and (2.22) will be used in this paper
to see the algebraic structure (the pure-state density operator) and
the $q$-deformation (the exponential density operator).
The coherent density operator was already used to interpret the physical
meanings of the exponential phase operator.
Then the relation (3.5) gives an algebraic relation depending on the forms of
the density operators. Generally we use the pure-state density operator only.

We briefly discuss the algebraic solution for  $B_\pm$ from (3.5).
We take the pure-state density operators for both $A_\pm$ and $B_\pm$.
The relation reduces to an algebraic difference equation in the variables
$n_a$  and $n_b$, formally written as
\be
{\cal A} (n_a+1) {\cal B}{} (n_b+1) \- {\cal A}{} (n_a){\cal B}{} (n_b)\=
{\cal D} (n_a, n_b).
\ee
Since the full structure of the $A_\pm$ is assumed known, we explicitly
know about the form of $\cal A$. The relation (3.6) becomes
a difference equation for $\cal B$. So we may find an algebraic solution
for ${\cal B}$, and sequentially obtain the form of $B_\pm$ with respect to
$n_b$. The $\cal B$ formally represents the square of magnitude of $B_\pm$.
Since we assume from the outset that $B_\pm$ changes a quantal number
of $N_b$, they are thus easily realized by the product of the square root
of $\cal B$ and the exponential phase operator.
The so-obtained solution of $B_\pm$ is a formal one, we should
change it into a normal form expressed by the boson generators $(b,\;\bd)$.
This can be done by using the definition of the exponential phase operator
(2.9) and (2.10).

But the situation is similar if we take an exponential density operator
for $A_\pm$ and a pure-state density operator for $B_\pm$ respectively.
We should choose the pure-state for $B_\pm$, for their full Hilbert space is
not assumed to be known.
The relation (3.5) under the resulting expectation value of $A_\pm$ gives
a new relation for $B_\pm$ depending on the parameter arising from the density operator.
In other words, the relation for $B_\pm$ looks like the deformed algebra
with respect to the algebraic solution.
In reality, the forms of $B_\pm$ as the solution should not
depend on the forms of the density operators, since the algebra of
$(D_\pm, D_0)$ is satisfied without the density operator.
The resulting algebra should not have the deformation parameter [19, 20].
But, in many cases, the relation takes the form of the well known
$q$-deformed algebra under redefinition of the parameter.
As a result, we can consider the process of taking the expectation value
as a tool to obtain a $q$-deformation physically.
Also this method explicitly and directly explains
why the mutual transformation in the equivalent class of the $q$-deformations
is possible, and why they take such forms.
The explanation is that they are related with each other
by the various transformation from the initial undeformed operators.

We now proceed to $q$-bosons.
The $q$-bosons have been defined in the different forms;
\bea
\tbm\tbp\- q^2\tbp\tbm &=& \hskip4cm\hbox{type \I} \hskip-5.0cm 1\\
                       &=& \hskip4cm\hbox{type \III} \hskip-5.05cm 1-q^2\\
                       &=& \hskip4cm\hbox{type \II} \hskip-5.1cm q^{-2N_b} \\
              &=& \hskip4cm\hbox{type \IV} \hskip-5.05cm (1-q^2)q^{-2N_b}
\eea
Here we impose the various structural types to the $q$-bosons for
our own convenience.
These $q$-bosons have been introduced by different authors [1 - 10].
Their mutual transformations can be found in [11 - 13].
We now reconstruct the above types of $q$-bosons
from the relation  (3.5).

First, we proceed to the $q$-boson of type \I by taking $A_\pm$
as the exponential phase operators and $D_0$ as the identity operator up to the
$a$-vacuum,
\be
A_+ \= e^\dagger_a ,\qq A_-\= e_a ,\qq D_0 \= 1. 
\ee
The relation (3.4) is rewritten as
$e_a e_a^\dagger B_- B_+\- e_a^\dagger e_a B_+ B_- \= 1.$
The algebraic solution of $B_\pm$ is approximately the bosonic generators,
\ie $B_+\=\bd,\;B_-\=b$.  We note that, to get an algebraic bosonic solution,
$D_0$ should be in the form $D_0\=1\+ N_b\mid 0_a\ra \la 0_a \mid$.
Thus we say that $D_0$ is the identity operator up to the $a$-vacuum.
Taking $D_0=1$ is {\em a kind of a constraint}.
For more detail treatment of the constraint see [19].

Then what is the meaning of $D_+\=e^\dagger_a B_+\:(D_-\=e_a B_-)$?
To see this, we recall that the exponential phase operator becomes the physical
phase operator in the large mean number limit ( the classical limit) of
the coherent state,
\ie $\la z |e_a|z\ra \= z/|z|\+\cdots\sim e^{i\f}$.
Since $B_\pm$ is equal to the boson generators, the related generators
$D_- =e_a B_-\;(D_+=e_a^\dagger B_+)$ are similar to a global $U(1)$ gauge
transformation of the normal boson generators in the classical limit
of the coherent state of the $a$-boson.
Thus we can say  that the operators $D_\pm$ are the operatal global $U(1)$
gauge transformation of the normal boson.
The $B_\pm$  is obtained from $D_\pm$ by the operatal $U(1)$ gauge fixing,
so $B_\pm$ is interpreted as the gauge fixed operator of $D_\pm$.
Taking the expectation value can be understood as gauge fixing.

The relation (3.5) under the condition (3.11), \ie taking the expectation
values for the exponential density operator of $A_\pm$ in (2.26) and (2.27),
gives
\be
\tbm \tbp \- q^2 \tbp\tbm \= 1.
\ee
Here, we put the tilde on the operator in order to represent the connection
between the deformed and the original operators.
As a result, the $\tilde{B}_\pm$ is the $q$-boson of type \I.
We note that the $q$-parameter arises from the exponential distribution of
the statistically mixed state.  We use the density operator;
\be
\rho_{en_b} \=(1-q^2) \sum_{n_a=0}^\infty q^{2 n_a}
\mid n_a,n_b \ra \la n_b,n_a\mid.
\ee
As a result, we can say that {\em the operatal $U(1)$ gauge transformed boson
is deformed into the $q$-boson of type \I}.

As a second example, we choose $A_\pm$ as the boson and $D_0$ as the identity
operator up to the $b$-vacuum in (3.5).
\be
A_+ \= \ad ,\qq A_-\= a ,\qq D_0 \=1 . 
\ee
This example is the reverse choice of the first example (3.11).
The relation (3.4) is rewritten as $a\ad B_- B_+\- \ad a B_+ B_-\=1$.
Note that $D_0\=1\+N_a\mid 0_b \ra \la 0_b \mid$.
Thus the state sum for $A_+=\ad\; (A_-=a)$ is just like that of the virtual
particle state,
and as a result finds the gauge information only.

The algebraic solution for $B_\pm$ is the exponential phase operator,
\be
B_- \= e_b,\qq B_+\=e_b^\dagger .
\ee
This solution is easily seen from the first example (3.11).

We take the expectation value with the density operator (3.13),
then consulting (2.24) and (2.25) gives the relation
\be
\tbm \tbp \- q^2 \tbp\tbm \= 1- q^2.
\ee
This relation shows the $q$-boson of type \III.
As a result, {\em the gauge transformation ,\ie the exponential phase operator,
is deformed into the $q$-boson of type \III}.

We now construct the $q$-boson of type \IV, and back to type \II.
We take $A_\pm$ as the boson and $D_0$ as the identity
except up to the $\a$-th state of the $a$-boson and the $b$-vacuum, \ie
\be
A_+ \= \ad ,\qq A_-\= a ,\qq D_0 \= \t(N_a \-\a).
\ee
The relation (3.4) is rewritten as
$a\ad B_- B_+ \- \ad a B_+ B_-\=\theta(N_a -\a)$.
In order to keep the algebraic structure, we require that $\a$
is independent of $n_a$,
\be
[N_a,\a]\= 0.
\ee
The corresponding algebra looks very strange in terms of $D_0$, since it is
equal to the identity depending on states of the algebra.
This special property will be treated further after constructing the $q$-boson
of type \IV.

The algebraic solution for $B_\pm$ of (3.4) is the exponential phase operator
satisfying (2.10). After taking the expectation value for the density operator
(3.13) and consulting (2.24), (2.25) and (2.28), the relation (3.4) gives
the $q$-boson,
\be
\tbm \tbp \- q^2 \tbp\tbm \= (1-q^2)q^{-2\a}.
\ee
This relation shows the $q$-boson of type \IV if we take
\be
\a \= N_b,
\ee
\be
\tbm\tbp \- q^2 \tbp\tbm \= (1-q^2) q^{-2N_b}.
\ee
This choice is acceptable, since $[N_a, N_b]\=0$.

From now on, we consider the meanings of this choice (3.17). Under the choice,
the relation (3.4) is written as
\be
a\ad B_-B_+ \- \ad a B_+ B_- \= \t(N_a \-\a).
\ee
We first change the nontrivial term $\t(N_a \- \a)$ into the identity
by a proper transformation.
It is just an operatal similarity (\ie adjoint) transformation of $x$ by
the $\a$-th power  of the exponential phase operator $e_a^\a  $
such as $e_a^{\dagger\a}\:x\: e_a^\a$.
We take the operatal similar transformation in which $D_0=\t(N_a-\a)$
changes into the identity operator,
\be
e_a^{\dagger\a} \t(N_a\-\a) e_a^\a \= \t(N_a) \=1.
\ee
We act an operatal similarity transformation on the boson algebra
$(a,\ad,N_a)$, and write them as $(a(\a),\ad(\a),N_a(\a))$ with
\be
a(\a) \=e_a^{\dagger\a} a e_a^\a,\qq \ad(\a) \= e^{\dagger\a}\ad e_a^\a.
\ee
These operators $(a(\a),\ad(\a)) $ form a {\em new} boson,
so we call it the {\em $\a$-adjoint boson} from now on.
The relation (3.22) changes into (3.14) under the $\a$-adjoint boson,
\ie
\be
a(\a)\ad(\a) B_-B_+ \- \ad(\a) a(\a) B_+ B_- \= 1.
\ee
The $q$-boson (3.21) is an operatal $\a$-adjoint transformed version of (3.16).

Let's consider more properties of the $\a$-adjoint boson.
The generators of the $\a$-adjoint  boson satisfy the boson algebra
with a different vacuum,
\be
[ a(\a), \ad(\a) ] \= \t(N_a -\a).
\ee
The vacuum of the $a(\a)$-boson is defined in the same way as the $a$-boson with
\be
a(\a) | 0_a(\a)  \ra \= 0.
\ee
This vacuum contains up to the $\a$-th number state of the normal boson,
\be
a(\a) | n_a \ra \= 0,\qq \hbox{if} \;\; n_a \leq \a.
\ee
Thus we can think of the $\a$-vacuum as the sum of the annihilated states,
\be
|0_a(\a) \ra \= \sum_{n_a=0}^\a c_{n_a} |n_a\ra,
\ee
where the $c_{n_a}$'s are constants that should satisfy the normalization
\be
\sum_{n_a=0}^\a |c_{n_a}|^2 \= 1,
\ee
since the vacuum is normalized to unity ($\la 0_a(\a) | 0_a(\a) \ra = 1$).
The vacuum of the $a(\a)$-boson forms the fundamental representation
of the global $SU(\a+1)$ symmetry.

We define the number operator $N_a(\a)$ of the $a(\a)$-boson by
\be
N_a(\a) \= \ad(\a) a(\a).
\ee
The number state $|n_a(\a) \ra$ of the $a(\a)$-boson is given by
\bea
N_a(\a) |n_a(\a) \ra &=& n_a |n_a(\a) \ra,\\[2mm]
|n_a(\a) \ra &=& \frac{1}{\sqrt{n_a(\a)!}} (\ad(\a))^{n_a} |0_a(\a)\ra
\= |n_a+\a\ra,
\qq n_a\geq 0.
\eea
We take the expectation value of the $a(\a)$-boson for the exponential density
operator (2.22).

Motivated from the above operatal adjoint transformation of the boson
generators, we act with the transformation on the exponential phase operator,
\be
e_a(\a)\= e_a^{\dagger\a} e_a e_a^\a, \qq
e^\dagger_a(\a)\= e_a^{\dagger\a} e^\dagger_a e_a^\a.
\ee
These exponential phase  operators act like the exponential phase  operator (2.10) for the
vacuum $|0_a(\a)\ra$,
\be
e_a(\a) e^\dagger_a(\a) \- e^\dagger_a(\a) e_a(\a) \=|0_a(\a)\ra\la 0_a(\a)|.
\ee
They also change the quantum number by one step and satisfy
the same form of commutation relations as in (2.11).

Finally we can construct the $q$-boson of type \II by using the $\a$-adjoint
exponential phase operator.
We take $A_\pm$ as the $\a$-adjoint exponential phase operators
and $D_0$  as the identity up to the vacuum $|0_a(\a)\ra$ in (3.4),
\be
A_+\=e_a^\dagger(\a),\qq A_-\=e_a(\a),\qq D_0 \= 1.
\ee
Then the algebraic solution for $B_\pm$ is the bosonic generators,
\ie $B_+ = \bd,\; B_- = b$.
Thus $D_\pm$ can be considered as the operatal global $U(1)$ gauge
transformation of the boson generators.
The relation (3.5), under the density operator (3.13),  gives
\be
\tbm \tbp \- q^2 \tbp\tbm \= q^{-2\a}.
\ee
This algebra is the $q$-boson of type \II if we choose $\a\=N_b$, with
\be
\tbm \tbp \- q^2 \tbp\tbm \= q^{-2N_b}.
\ee
As a result, the $q$-boson of type \II is related to  the gauge transformation
by the adjoint transformation of the exponential phase operator.

\section{$SU_q(N)$-covariant bosons and $SU(N)$ symmetry}

We now proceed to extend the system (3.1) to an $N$-component system
$(D_{+i}, D_{-i}, D_{0ij})$ which satisfies the commutation relations
\be
[D_{-i}, D_{+j}] \= D_{0ij}  \label{4.1}
\ee 
in the formal sense.
All other commutators are related to the ladder structure of the algebra,
and we will fix them in the middle of the construction of $q$-bosons.

Similar to the previous one component case,
we introduce two sets of $N$ component independent bosons
$\{ (a_i, a^\dagger_i, N_{ai}), \;(b_j, b^\dagger_j, N_{bj})\}$ such that
\bea
[a_i , a^\dagger_j ] \= \delta_{ij}\= [ b_i , b^\dagger_j ],
\qq i,j=1,2,\cdots,N,\nn\\[1mm]
N_{ai} \= a_i^\dagger a_i, \qq N_{bi} \= b^\dagger_i b_i, \qq i=1,\cdots,N.
\eea 
The Hilbert space, characterized by the number operators $(N_{ai}, N_{bi})$,
is spanned by pure states  $\mid n_{ai},\cdots,n_{bN} \rangle $.

We assume that $D_{\pm i}$ are decomposed into the $a_i$- and $b_i$-contributions,
\be
D_{-i} \= A_{-i} B_{-i}, \qq D_{+i} \= A_{+i} B_{+i}.
\label{4.4}
\ee  
The mutually independent $A_{\pm i}$ and $B_{\pm i}$ are fixed to
the one step operator,
\bea
[ N_{ai} , A_{\pm j} ] \= \pm \d_{ij} A_{\pm i},
 & [ N_{bi} , B_{\pm j} ] \= \pm \d_{ij} B_{\pm i} , \nn\\[1mm]
[A_{\pm i} , A_{\pm j}]\=0, & [A_{\pm i}, B_{\pm j} ]\=0\=[A_{\pm i}, B_{\mp j}],
\qq i\neq j.
\label{4.5}
\eea 
For simplicity, the choice of $A_{\pm i}$ is restricted to a normal boson and
an exponential phase operator.
Also, $D_{0ij}$ is fixed to the identity operator and the step function
to find the known $q$-boson.
After substituting the product form of $D_{\pm i}$ into the commutator
\mr{4.1}, we rewrite the relation as
\be
A_{-i}A_{+j} B_{-i} B_{+j} \- A_{+j}A_{-i} B_{+j} B_{-i} \= D_{0ij}.
\label{4.6}
\ee 
Here $D_{0ij}$ is a function of the number operators $(N_{ai}, N_{bi})$.
We take an expectation value of the relation \mr{4.6} for an exponential
density operator for the $a_i$-boson and a pure-state density operator for the
$B_{\pm i}$ state, with
\be
\rho\=\sum_{n_{ai}=0}^{\infty}\Big(\prod_{i=1}^N (1-q_i^2)\: q_i^{2n_{ai}}\Big)
      \mid n_{ai}, n_{bi} \ra \la n_{ai} , n_{bi} \mid.
\label{4.7}
\ee 
Then the result will be a $q$-deformation of the algebra
of $B_{\pm i}$ as we have done in the previous section.
The $q$-parameters are given by
\be
q^2_i \= e^{-\e_i/k_BT}, \qq i=1,2,\cdots,N.
\label{4.8}
\ee 
Note that these $q$-parameters differ with the quantal energies.

We first take the exponential phase operator for $A_{\pm i}$
and the identity operator for $D_{0ij}$ in (4.5) as
\be
A_{-i}\=e_{ai}(\a_i),\; A_{+i}\=e^\dagger_{ai}(\a_i),\qq D_{0ij} \= \d_{ij}.
\label{4.9}
\ee 
Here, $e_{ai}(\a_i)$ is the similarity (adjoint) transformation of
the exponential phase operator,
\be
e_{ai}(\a_i)\= e_{ai}^{\dagger\a_i} \:e_{ai}\: e_{ai}^{\a_i}.
\ee 
In order for the $A_{\pm i}$ to be indepedent on each other and the $B_{\pm j}$,
$\a_i$ should commute with the number operators $N_{ai}$  and $N_{bj}$.
Furthermore, we take them as numbers.
The relation (4.5) is rewritten as
$e_{ai}(\a_i)e^\dagger_{aj}(\a_j) B_{-i}B_{+j}\-e^\dagger_{aj}(\a_j) e_{ai}(\a_i)
B_{+j}B_{-i}\=\d_{ij}$.

The algebraic solutions for $B_{\pm i}$ in (4.5) are bosons,
\be
B_{-i} \= b_i, \qq B_{+i} \= b_i^\dagger.
\ee 
After taking the expectation value \mr{4.6} with respect to (4.6),
the nontrivial relations are seen in the same species, since the density operator
(4.6) is diagonal, with
\be
\tilde{B}_{-i} \tilde{B}_{+i} \- q_i^2 \tilde{B}_{+i}\tilde{B}_{-i}
\=q_i^{\a_i}.
\label{4.12}
\ee 
These are relations of type \II for each species.
We then encounter a problem to determine the other relations between
the different species consistently.
To treat this problem properly, we should take a $q$-differential calculus
and a Yang-Baxter equation into consideration [29,30].
But we leave further details to a future paper,
treat only the independent cases, and
restrict our attention here to finding the hidden symmetries of the system.

We assume that all the $\a_i$ are equal to zero,
\be
\a_i\=0.
\label{4.13}
\ee 
The different species of $B_{\pm i}$ are mutually independent and commuting in
the relations \mr{4.12}.
Thus we construct the $N$-component independent $q$-bosons of different
$q$-deformation parameters,
\bea
\tbmi \tbpj \- q_j^{2\d_{ij}} \tbpj\tbmi &=& \d_{ij}, \\[1mm]
[N_{bi}, \tilde{B}_{\pm j} ] &=& \pm \tilde{B}_{\pm j} \d_{ij}  \nn
\eea 
These independent $q$-bosons are simple extensions of the $q$-boson of type
\II into  $N$ components.

We construct the $su_q(N)$ algebra by using the Schwinger method, and find
the hidden symmetry taking place in the process of constructing the algebra.
The $su_q(N)$ algebra is given by
\bea
[H_i,H_j]\=0, & [E_i,F_j]\=\delta_{ij} \: [H_i], \nn\\[1mm]
[H_i,E_j]\=\hbox{A}_{ij} E_j, & [H_i,F_j]\=-\hbox{A}_{ij} F_j,
\eea 
where $[x] \equiv (q^x - q^{-x}) / (q - q^{-1})$, and
$\hbox{A}_{ij}\=2\delta_{ij}\-\delta_{i j+1}\-\delta_{ij-1}$ is
the element of the $su(N)$ Cartan matrix.

This algebra has only one deformation parameter.
In order to construct $su_q(N)$ by combining the two independent
$q$-bosons of (4.13) in the product form,
we should require all $q_j$-values to be equal, with
\be
q_1^2\=q_2^2\=\cdots\=q_N^2 \equiv q^2.
\ee 
The same $q$-parameter means that the quantal energies of the hamiltonian
in the density operator \mr{4.7} should be equal,
that is $\e_1=\e_2=\cdots=\e_N=\e$. Thus the hamiltonian becomes
\be
H\=\e\sum a_i^\dagger a_i.
\ee
The hamiltonian in the density operator (4.6) has a global $SU(N)$ symmetry,
\bea
a_i^\prime\=\sum_j U_{ij} a_j, \qq
a_i^{\prime\dagger}\=\sum_j U^\dagger_{ij} a_j^{\prime\dagger},\nn\\[1mm]
\sum_k U_{ik} U^\dagger_{kj} \=\sum_k U^\dagger_{ik} U_{kj}\=\d_{ij}.
\eea
Under this symmetry, the $q$-boson (4.13) is changed into a $q$-parameter algebra;
\bea
\tbmi \tbpj \- q^{2\d_{ij}} \tbpj\tbmi &=& \d_{ij}, \\[1mm]
[N_{bi}, \tilde{B}_{\pm j} ] &=& \pm \tilde{B}_{\pm j} \d_{ij}.  \nn
\eea
This algebra has the global $SU(N)$ transformation (4.17) as its hidden symmetry.
Also the Chevalley basis of $su_q(N)$  are given, in [4], by
\be
H_i\=N_i-N_{i+1}, \; E_i \= \tbpi\tilde{B}_{-(i+1)}, \;
F_i\= \tilde{B}_{+(i+1)}\tbmi, \; i = 1,\cdots, N-1.
\ee 
As a result, we find the relation between the global $SU(N)$ symmetry and
the Schwinger realization of  $su_q(N)$.

As a second example, we consider the $SU_q(N)$-covariant bosons [5], where
\bea
\tbmi \tbpi \- q^2 \tbpi \tbmi \= q^{2\sum_{j=1}^{i-1} N_j}, \nn\\[1mm]
\tbmi \tbmj \= q \tbmj \tbmi, \qq i < j, \\[1mm]
\tbmi \tbpj \= q \tbpj \tbmi \qq i \neq j.\nn
\eea
The algebra can be rewritten by the $SU(N)$ $R$-matrix,
\be
R \= q \sum_i e_{ii}\otimes e_{ii}\+ \sum_{i\ne j}e_{ii}\otimes e_{jj}
\+ (q-q^{-1}) \sum_{i<j}e_{ij} \otimes e_{ji},
\ee
where $e_{ij}$ is the $N\times N$ matrix with entry $1$ at position $(i,j)$
and $0$ at elsewhere.
Using the notation $R\= R_{ij,kl} e_{ik}\otimes e_{jl}$, after a minor
modification, the relations (4.20) are rewritten as
\begin{eqnarray}
\tbmi \tbmj &=& q^{-1} R_{ij,kl} \tilde{B}_{-l} \tilde{B}_{-k} , \nn \\
\tbpi \tbpj &=& q^{-1} R_{lk,ij} \tilde{B}_{+k} \tilde{B}_{+l},  \\
\tbmi \tbpj &=& \delta_{ij} \+ q R_{ki,jl} \tilde{B}_{+k}
\tilde{B}_{-l}. \nn
\end{eqnarray}

We now derive the algebra (4.20).
We choose new forms of $A_\pm$ and $D_{0ij}$ in (4.5) as follows,
\be
A_{-i}\=e_{ai},\qq A_{+i}\=e^\dagger_{ai},\qq
D_{0ij} \= \theta(N_{ai}\- \sum_{k=1}^{i-1} N_{bk}) \d_{ij}.
\ee
The relation \mr{4.6} can be written as
\be
e_{ai}e_{aj}^\dagger B_{-i}B_{+j} \- e_{aj}^\dagger e_{ai} B_{+j}B_{-i}
\= \delta_{ij} \: \theta(N_{ai} \- \sum_{k=1}^{i-1} N_{bk}).
\ee
We take the expectation value for the density operator (4.6) and obtain the
nontrivial relation only for $i=j$.
We also require that the $a_i$-bosons are all equivalent, \ie their Hilbert
space has the global $SU(N)$ symmetry as in the independent $q$-bosons.
The nontrivial part under the expectation values of (4.24) are given by
\be
\tbmi \tbpi \- q^2 \tbpi \tbmi \= q^{2\sum_{j=1}^{i-1} N_j}.
\ee
The different $q$-bosons may be mutually commuting if we ignore
the tower dependence on the number operators (4.25).
But the right hand side of (4.25) prevents them from being independent
$q$-bosons.
We introduce independent (commuting) operators which take the right hand
side of (4.25) into a constant, with
\be
\hat{B}_{\pm i} \= q^{\sum_{k<i}N_k} \tilde{B}_{\pm i}.
\ee
Then we obtain the algebra (4.18). We can thus treat them as independent
in order to fix the relation between different species,
\bea
\hat{B}_{\pm i} \hat{B}_{\pm j} &=& \hat{B}_{\pm j} \hat{B}_{\pm i}, \nn \\
\hat{B}_{\pm i} \hat{B}_{\mp j} &-& q^{2\d_{ij}} \hat{B}_{\mp j}
\hat{B}_{\pm i} \= \d_{ij}.
\eea
These relations can be rewritten in terms of the original operators $\tilde{B}_{\pm i}$.
We then obtain the $q$-commuting property
\bea
\tbmi \tbmj \= q \tbmj \tbmi,\qq \tbpi\tbpj\= q\tbpj\tbpi\qq  i < j, \nn\\[1mm]
\tbmi \tbpj \= q \tbpj \tbmi \qq i \neq j.
\eea
After combining (4.25) and (4.28), we obtain the fundamental representation
of $SU_q(N)$. In reality, we needed to build in the commutatitivity
of the transformed operators of the different species by hand.
A more systematic and general approach requires the $q$-differential calculus [29, 30].

We now consider the meaning of the choice of (4.24) as done in the previous section.
Introduce the operatal similar transformation
\be
e_{ai}(\a_i)\= e_{ai}^{\dagger\a_i}\: e_{ai} \: e_{ai}^{\a_i},\qq
e_{ai}^\dagger(\a_i)\= e_{ai}^{\dagger\a_i}\: e_{ai}^\dagger \: e_{ai}^{\a_i},
\ee
where
\be
\a_i\=\sum_{k=1}^{i-1} N_{bk}.
\ee
Then the system  (4.24) is changed into
\be
e_{ai}(\a_i) e^\dagger_{aj}(\a_j) \: B_{-i}B^\dagger_{+j}
       \- e^\dagger_{aj}(\a_j)e_{ai}(\a_i) \: B_{+j} B_{-i} \= \delta_{ij}
\ee
This relation is equivalent to (4.8) except for the $\a$-adjoint exponential phase
operators.

As a result, the $SU_q(N)$-covariant bosons are related with the adjoint transformation
of the exponential phase operator.

\section{Discussions and Conclusions}

We see that the averaging method is a very convenient way to find the difference and
similarity between the different forms of $q$-bosons.
But to obtain our results, we restricted the product of two commuting generators
to satisfy the algebra which is different from the product algebra up to
their vacuum structure.
This is a kind of constraint in which the $q$-boson is based, so we need to
consider it further on the basis of the Hopf algebra [16, 17].
Although we were able to extend our ideas to the multicomponent system
and find the hidden symmetry, it still remains to fully consider it using
the machinery of the $q$-calculus [29, 30].

There are many different $q$-bosons, not treated here.
However, they can be included in our method by relaxing and modifying
the conditions given here.
As a simple example, we can construct, by using a similar method,
the $q$-fermion [6,9],
\[
\tilde{C}\tilde{C}^\dagger\+ q^{2} \tilde{C}^\dagger \tilde{C}\= 1
\]
More interesting things appear when the density operator is fermionic
and this system is related to extending the real $q$-parameter
into a complex variable.

Various transformations using the exponential phase operator show some
interesting properties, since it is a $U(1)$ phase in the classical limit.
We can thus obtain an operatal global $U(1)$ gauge transformation
( see the $q$-boson of type \I) by taking products with an exponential phase operator.
We can also define an operatal adjoint (similarity) transformation.
But such transformations are not unitary,
since the exponential phase operator is not unitary.
In particular, the adjoint transformation of a boson algebra (3.26) has the shifted
vacuum of a global $SU(\a)$ symmetry (3.29).

In order to generate a non-abelian gauge transformation, we
need a non-abelian phase operator that is still not defined.
More generally, it will also be worthwhile to extend the above concepts into
the field theory [26] and the matrix model [23 - 25].

\vskip1cm\noindent
{\bf \large Acknowledgment}
\vskip.5cm

\noindent
The author thanks Murray Batchelor for his hospitality and help with
the manuscript.

\end{document}